\documentclass[graybox]{svmult}

\usepackage{type1cm}        
\usepackage{makeidx}         
\usepackage{graphicx}        
\usepackage{multicol}        
\usepackage[bottom]{footmisc}

\usepackage{newtxtext}       %
\usepackage{newtxmath}       
\usepackage{multirow}
\usepackage{diagbox}

\usepackage{xcolor}
\usepackage[colorlinks]{hyperref}

\usepackage{booktabs}

\usepackage{amsmath,amsfonts}%
\usepackage{amsthm}%
\usepackage{mathrsfs}%
\usepackage[title]{appendix}%
\usepackage{textcomp}%
\usepackage{manyfoot}%
\usepackage{booktabs}%
\usepackage{algorithm}%
\usepackage{algorithmicx}%
\usepackage{algpseudocode}%
\usepackage{listings}%
\usepackage{bm}%
\usepackage{enumitem}
\usepackage{caption}

\theoremstyle{thmstyleone}%

\definecolor{darkgreen}{rgb}{0.0,0.4,0.0}
\definecolor{darkred}{rgb}{0.6,0.0,0.0}
\definecolor{darkblue}{rgb}{0.0,0.0,0.5}
\definecolor{gray}{rgb}{0.5,0.5,0.5}
\definecolor{cyan}{rgb}{0.0,1.0,1.0}
\definecolor{darkcyan}{rgb}{0.0,0.5,0.5}
\definecolor{darkorange}{rgb}{0.8,0.4,0.0}
\definecolor{darkmargenta}{rgb}{0.5,0.0,0.5}
\definecolor{black}{rgb}{0.0,0.0,0.0}

\hypersetup{citecolor=blue, linkcolor=red, anchorcolor=darkblue,
filecolor=darkblue, urlcolor=darkblue}
\hypersetup{pdfauthor={Nicolae Suciu},pdftitle={Suciu}}

\usepackage{float} 
\usepackage{bm}
\usepackage{subcaption}

\makeindex             

\def \e  {\varepsilon}

\def \z  {{\boldsymbol z}}

\def \div {\nabla\cdot}

\def \bmu {\bm{u}}
\def \bmv {\bm{v}}

\def\bmf{\boldsymbol f}

\def\bmu{\boldsymbol u}
\def\bmv{\boldsymbol v}

\begin{document}

\title*{Fixed stress splitting approach for contact problems in a porous medium}
\author{Tameem Almani\orcidID{0000-0002-5603-2533} and\\ Kundan Kumar\orcidID{0000-0002-3784-4819}}
 \institute{Tameem Almani,  \at Saudi Arabian Oil Company (Saudi Aramco), Saudi Arabia, \email{tameem.almani@aramco.com}
\and Kundan Kumar \at Center for Modeling of Coupled Subsurface Dynamics, University of Bergen, Bergen, Norway, \email{kundan.kumar@uib.no} }
\maketitle


\abstract{We consider a poromechanics model including frictionless contact mechanics. The resulting model is the Biot equation with contact boundary conditions leading to a variational inequality modelling mechanical deformations coupled to a linear parabolic flow equation. We propose a fully discrete iterative scheme for solving this model, extending the fixed-stress splitting scheme. We use finite elements in space and a backward Euler discretization in time. We show that the fixed stress split scheme is a contraction. }

\section{Introduction}
\label{sec:intro}
In this article, we present a splitting method for solving a contact mechanics problem for a porous medium, with a particular focus on the role of fluid. Our motivation to study this problem stems from the role of contact mechanics in a fractured subsurface medium. Assuming the fracture pressure to be known and symmetry of the geometry, we can reduce the problem of contact mechanics of the fracture to that of poroelasticity with contact mechanics. Mathematically, the model takes the form of the Biot equations \cite{Biot1955} with contact conditions \cite{sofonea2012mathematical} at a part of the boundary. We refer to the recent works of \cite{banz2023contact, hosseinkhan2021biot, hosseinkhan2023semilinear} where a detailed mathematical and numerical analysis of the contact problem for the Biot model have been carried out. 

Coupling of flow and mechanical deformations in the subsurface is relevant for a wide range of applications, including hydrocarbon recovery, geothermal energy production, and geological hazard assessment. A key challenge is in the understanding of mechanics of the faults and fractures especially in the presence of fluids. The subsurface geoengineering processes involve injection or extraction of fluids leading to a build up of pressure and possible imbalance of in-situ stress conditions. This may result in seismic events. Understanding the mechanics of fractures and their interaction with fluids is therefore a topic of great importance in geosciences. This also stems from the fact that the faults and fractures are mechanically the most vulnerable regions for mechanical failures. At the same time, they have a strong influence on the flow profiles. Due to the heterogeneities and roughness of the fractures and fault surface, their presence can lead to complex mechanical behavior. In the presence of fluid, the mechanical behaviour gets further complicated, as the fluid pressure reduces the effective compressive normal stresses making it easier to slip. Moreover, the fluid behaviour can be quite complicated as they can act as a lubricant and reduce the frictional forces that resist sliding or shear. The study of contact mechanics of fractured porous medium coupled to the role of fluid is of high relevance.\par

Contact problems in a poroelastic medium leads to a coupled system of flow and mechanical deformations. Iterative methods are suitable for such problems as they decompose the problem into its individual constituents of flow and mechanical problems. The individual equations can be accordingly solved using well-developed solvers and pre-conditioners. Moreover, the splitting methods allow flexibility in terms of time stepping \cite{Almanietal2016} for each of these equations. However, design of such iterative schemes require careful considerations as poorly designed schemes may lead to unstable schemes. Two main iterative coupling schemes to solve the flow problem coupled with geomechanics  are the fixed-stress split and the undrained split schemes \cite{SettariMourits1998,Kimetal2009, GiraultKumarWheeler2016}.   In the linear Biot case, both schemes were shown to be convergent in 
\cite{MikelicWheeler2012,  AlmaniKumarWheeler2017} and \cite{TAKA}. Here we focus on extending the fixed stress split scheme to the current setting. Our main contribution is in 
the convergence analysis of the fixed stress split scheme. Our convergence analysis proof uses the approaches in \cite{MikelicWheeler2012, AlmaniKumarWheeler2017} that have been developed for the linear Biot equations.

\section{Model equations}
 Let $\Omega \subset \mathbb{R}^d, d = 2,3$ be a bounded domain with smooth boundary $\Gamma$. We assume 
$\Gamma = \Gamma_1 {\displaystyle \cup } \Gamma_2 {\displaystyle \cup } \Gamma_3$, with $\Gamma_1, \Gamma_2, \Gamma_3$ being pairwise disjoint sets. In addition, $\text{meas}(\Gamma_1) >0.$ Let $\mathbb{S}^d$ denote the space of second order symmetric tensors on $\mathbb{R}^d$. The domain $\Omega$ is a porous medium. We consider the following problem where we consider the quasi static Biot model in $\Omega$ with contact conditions at part of the boundary $\Gamma_3$ including the frictionless case with the Signorini condition. We refer to \cite{banz2023contact} for further details on the model and its mathematical and numerical analysis. \par

{\textbf{Problem}}: (Frictionless contact with the Signorini condition in a form with a gap function in a porous medium.)
Let us first consider the equation for mechanical displacement. Find a displacement field $\bmu: \Omega \mapsto \mathbb{R}^d$ and a stress field $\boldsymbol \sigma: \Omega \mapsto \mathbb{S}^d$ such that
\begin{align}
   \label{eq:mechanics1} \boldsymbol \sigma^{\rm{por}}(\bmu, p) &= \boldsymbol \sigma(\bmu) - \alpha p \boldsymbol I,\\
    \label{eq:mechanics2} \boldsymbol \sigma(\bmu) &= \lambda (\nabla \cdot \bmu) I + 2G \varepsilon(\bmu),\\
    \label{eq:mechanics3} \boldsymbol \varepsilon(\bmu) &= \dfrac{1}{2} \left (\nabla \bmu + \nabla \bmu^T \right),\\
    \label{eq:mechanics4} \div \boldsymbol \sigma^{\rm{por}} + \bmf_0 &= 0.
\end{align}

The first equation \eqref{eq:mechanics1} is the definition of effective poroelastic stress field which incorporates the modification due to fluid pressure with $\alpha$ being the Biot coefficient. The second equation \eqref{eq:mechanics2} is  Hooke's law that specifies the linear constitutive relationship between stress and strain. Next equation \eqref{eq:mechanics3} is the definition of linear strain and the last equation \eqref{eq:mechanics4} is the conservation of linear momentum.  The fluid pressure $p$ is also an unknown here. We will supply the relevant equation for fluid pressure later. \par
 Next, we specify the boundary conditions. 
\begin{align}
  \label{eq:contact1}  \bmu = 0 \; \; \text{ on } \; \; \Gamma_1, \; \; \boldsymbol \sigma \boldsymbol \nu = \boldsymbol f_2 \quad \text{ on } \quad \Gamma_2,\\
  \label{eq:contact2}  u_\nu \leq g_a, \; \sigma_\nu \leq 0, \; \sigma_\nu (u_\nu - g_a) = 0 \quad \text{ on } \quad \Gamma_3, \\
   \label{eq:contact3} \boldsymbol \sigma_\tau = 0 \quad \text{ on } \quad \Gamma_3.
\end{align}
Here, we assume $g_a \geq 0$ is a given scalar (denotes the gap between the  domain and the foundation) function; $u_\nu = \bmu \cdot \boldsymbol \nu, \sigma_\nu = (\boldsymbol \sigma \boldsymbol \nu) \cdot \boldsymbol \nu$, $\boldsymbol \sigma_\tau = \boldsymbol \sigma \boldsymbol \nu - \sigma_\nu \boldsymbol \nu$.   The two boundary conditions in \eqref{eq:contact1} are quite standard; a homogeneous Dirichlet boundary condition is specified on $\Gamma_1$ and normal component of stress (that is, external force) is specified on the boundary $\Gamma_2$. It is the contact conditions on $\Gamma_3$ given in \eqref{eq:contact2} and \eqref{eq:contact3} which are of particular interest. \eqref{eq:contact3} denotes the frictionless case with the tangential component of the stress being zero. This is a simplification of the more realistic conditions involving friction which will be addressed elsewhere. The contact condition \eqref{eq:contact2} represents the so-called Signorini contact condition with a gap function $g_a$. The first condition restricts the amount of normal component of displacement to the initial gap $g_a$ implying contact with a rigid foundation. The second condition suggests that the normal stress is compressive. The third condition is a complementarity condition which says that if the normal displacement is less than the initial gap, the normal stress vanishes. On the other hand, when the gap is zero, there is no restriction on the amount of normal stress. This ensures that the contact never penetrates the surface of rigid foundation. For more details on the model and the mathematical foundations, we refer to the textbook \cite{sofonea2012mathematical}.\par

The flow equation consists of the mass balance equation: 
\begin{equation}
 \frac{\partial}{\partial t}\Big(\big(\frac{1}{M} + c_f \varphi_0\big) p  +  \alpha \nabla \cdot \bmu \Big) + \nabla \cdot \z  = q.
\label{eq:masscons1}
\end{equation}
Here, $p$ is the pressure and  $\z$ is the fluid flux.  The coefficient in the first term is total compressibility and accounts for the solid matrix compressibility given by $1/M$ and the fluid compressibility $c_f \varphi_0$. This is a linearized version and together with the linear mechanics equation \eqref{eq:mechanics1} - \eqref{eq:mechanics4} forms the quasi-static Biot model.  The flux $\z$ is described by the Darcy law
\begin{align}
 \label{eq:darcy}   \z = - \frac{1}{\mu_f} \boldsymbol K \big(\nabla\,p - \rho_{f,r} g \nabla\, \eta\big).
\end{align}
Let us define, the space for $\bmu$ and pressure $p$: 
\begin{align*}
    V_u &:= \left \{\bmv \in H^1(\Omega)^d: \bmv = 0 \text{ on } \Gamma_1, v_\nu \leq g_a \text{   on } \Gamma_3 \right \}, \\
V_p &:= \left \{p \in H^1(\Omega): p = 0 \text{ a. e.  on } \Gamma  \right \}. 
\end{align*}

\subsection{Mixed variational formulation}
\label{sec:Varform}
We will use a mixed formulation for the flow and conformal Galerkin formulation for the mechanics equation. The mixed method defines flux as a separate unknown and rewrites the flow equation as a system of first order equations. For the fully discrete formulation (discrete in time and space), we use a mixed finite element method for space discretization and a backward-Euler time discretization. Let $\mathfrak{T}_h$ denote a regular family of conforming triangular elements of the domain of interest, $\overline{\Omega}$. Using the lowest order Raviart-Thomas  (RT) spaces, we have the following discrete spaces (${V}_h $ for discrete displacements, $Q_h$ for discrete pressures, and ${Z}_h$ for discrete velocities (fluxes)):
\vspace{-0.2cm}
\begin{align}
& \hspace{1cm}{V}_h = \{\bmv_h \in V_u\,;\, {\forall} T \in  \mathfrak{T}_h, {{\bmv}_h}_{|T} \in {{\mathbb{P}}_1}^d\}, \\
& \hspace{1cm}Q_h = \{p_h \in L^2(\Omega) \,;\, {\forall} T \in  \mathfrak{T}_h, {p_h}_{|T} \in {{\mathbb{P}}_0}, p_h = 0 \,\mbox{on}\ \, \Gamma \}, \\
& \hspace{1cm} {Z}_h = \{{\boldsymbol q}_h \in {H(\text{div}; \Omega)}^d\,; {\forall} T \in  \mathfrak{T}_h, {{{\boldsymbol q}}_h}_{|T} \in {{\mathbb{P}}_1}^d, \, {{\boldsymbol q}}_h\cdot \nu = 0\,\ \mbox{on}\ \Gamma  \}.
\end{align}
We write down a fixed stress splitting scheme for the fully discrete Biot model with frictionless contact. Let us first derive a weak variational formulation for the mechanics equation \eqref{eq:mechanics1} - \eqref{eq:mechanics4} with the boundary conditions \eqref{eq:contact1} - \eqref{eq:contact3}. \par 
Following the standard procedure, multiplying with any element $\bmv - \bmu \in V_u$, our weak formulation finds $\bmu \in V_u$ such that
\begin{align*}
    2G \left ( \boldsymbol \e(\bmu), \boldsymbol \e(\bmv) - \boldsymbol \e(\bmu) \right )_\Omega + \lambda \left ( \div \bmu, \div \bmv - \div \bmu \right ) - \alpha \left(p, \div \bmv - \div \bmu \right) \\
    = (\sigma_\nu, v_\nu - u_\nu)_{\Gamma_3} + (\boldsymbol f_0, \bmv - \bmu) + (\boldsymbol f_2, \bmv - \bmu)_{\Gamma_2}+ (\boldsymbol \sigma_\tau, \bmv_\tau - \bmu_\tau)_{\Gamma_3} 
\end{align*}
for all $\bmv \in V_u$.
Using frictionless case, $\boldsymbol \sigma_\tau = 0$ on $\Gamma_3$, the last term drops out. 
We take a look at the contact mechanics term: 
\[(\sigma_\nu, v_\nu - u_\nu)_{\Gamma_3} = (\sigma_\nu, (v_\nu - g_a) - (u_\nu - g_a))_{\Gamma_3} = (\sigma_\nu, (v_\nu - g_a) )_{\Gamma_3}\]
by using the condition $\sigma_\nu (u_\nu - g_a) = 0$. With the other two conditions, $v_v - g_a \leq 0, \sigma_\nu \leq 0$, we deduce $(\sigma_\nu, (v_\nu - g_a) )_{\Gamma_3} \geq 0.$ The weak formulation for our problem thus takes the form: find $\bmu \in V_u$ such that 
\begin{align*}
    2G \left ( \boldsymbol \e(\bmu), \boldsymbol \e(\bmv) - \boldsymbol \e(\bmu) \right )_\Omega + \lambda \left ( \div \bmu, \div \bmv - \div \bmu \right ) - \alpha \left(p, \div \bmv - \div \bmu \right) \\
    \geq (\boldsymbol f_0, \bmv - \bmu) + (\boldsymbol f_2, \bmv - \bmu)_{\Gamma_2}
\end{align*}
for all $\bmv \in V_u$.

We adopt the usual notation for the discrete approximations. For the discrete case, with $k$ denoting the time index, $\Delta t$ the time step, the weak formulation for the mechanics equation takes the form of a variational inequality.\par
{\textbf{Discrete variational inequality for mechanics equation:}}
Find $\bmu_h^k \in V_h$ such that for all $v_h \in V_h$
\begin{align*}
2G \left ( \boldsymbol \e(\bmu_h^k),  \boldsymbol \e(\bmv_h) - \boldsymbol \e(\bmu_h^k )   \right ) + \left ( \lambda \div \bmu_h^k,   \div \bmv_h - \div \bmu_h^k  \right ) \\
- \alpha \left ( p_h^k,   \div \bmv_h - \div \bmu_h^k \right ) \geq  \left (\boldsymbol f_0,  \bmv_h - \bmu_h^k \right) +  \left (\boldsymbol f_2,  \bmv_h - \bmu_h^k \right)_{\Gamma_2}.
\end{align*}
This is coupled to the weak formulation for the flow problem. \par

({\bf Flow equation in discrete form})  Find  $p_{h}^{k} \in Q_h$ and ${\z}_{h}^{k} \in {Z}_h$ such that 
 
\begin{align}
\nonumber \forall {\theta}_h \in Q_h\,,\, &\frac{1}{\Delta t} \Big(\big(\frac{1}{M} + c_f \varphi_0  \big) \Big( p_{h}^{k}-p_{h}^{k-1} \Big), {\theta}_h \Big) + \frac{1}{\mu_f} \big(\nabla \cdot {\z}_{h}^{k}, {\theta}_h) =  \\
\label{eq:flowp_sr_fs_sr} & - \frac{\alpha}{ \Delta t}  \Big(\nabla \cdot  \Big( {\bmu}_{h}^{k}-{\bmu}_{h}^{k-1} \Big),  {\theta}_h \Big) + \Big({q}_h, {\theta}_h \Big), \\
\label{eq:flowf_sr_fs_sr} \forall {\boldsymbol q}_h \in {Z}_h\,,\, &\Big( {\boldsymbol K}^{-1} {\z}_{h}^{k}, {\boldsymbol q}_h \Big) =  \Big( p_{h}^{k}, \nabla \cdot {\boldsymbol q}_h \Big) + \Big(\rho_{f,r} g \nabla\,\eta, {\boldsymbol q}_h \Big).
\end{align}

We also need to specify the initial condition for the $1^{st}$ time step: $p_h^0 = p_0,\bmu_h^0 = \bmu_0. $
We assume the compatibility of the initial data, see \cite{Showalter2000}.   \par
{\textbf{Notation:}} In what follows, $n$ denotes the iteration index whereas $k$ denotes the time stepping index. In the scheme and the subsequent proof of contraction, $k$ is assumed to be fixed whereas the iteration index $n$ varies. As $n$ goes to infinity, the scheme converges to the fully implicit scheme.
\vspace{-0.8cm}
\subsection{Fixed stress split for frictionless contact problem}
 The scheme consists of two steps, flow solve followed by mechanics solve. In the flow solve, we add a regularization term and show that this is sufficient to ensure contraction.

Step (a):  Flow step: Find $p_{h}^{n+1,k} \in Q_h$, ${\z}_{h}^{n+1,k} \in {Z}_h$ such that:
\begin{align}
\nonumber \forall {\theta}_h \in Q_h\,,\, \Big(\big(\frac{1}{M} + c_f \varphi_0 + \frac{\alpha^2}{\lambda} \big) \big(\frac{p_{h}^{n+1,k} - p_{h}^{k-1}}{\Delta t}\big), {\theta}_h \Big)  
+ \frac{1}{\mu_f}\big(\nabla \cdot {\z}_{h}^{n+1,k},{\theta}_h \big) = \\
\label{eq:nf_eq4_fs_sr} \Big(-\frac{\alpha}{\lambda}\big(- \alpha\big(\frac{\textcolor{black}{p_{h}^{n,k}} - p_{h}^{k-1}}{\Delta t}\big)+\lambda \nabla \cdot \big(\frac{{\bmu}_{h}^{n,k} - {\bmu}_{h}^{k-1}}{\Delta t}\big)\big), {\theta}_h \Big)+ \big(q_h, {\theta}_h \big)
\end{align}
\begin{equation}
\forall {\boldsymbol q}_h \in {Z}_h\,,\, \big(\boldsymbol K^{-1}{\z}_{h}^{n+1,k},{\boldsymbol q}_h \big) = (p_{h}^{n+1,k}, \nabla \cdot {\boldsymbol q}_h) +\big (\nabla(\rho_{f,r} g \eta),{\boldsymbol q}_h)
\label{eq:nf_eq5_fs_sr}
\end{equation}

Step (b): Mechanics step. Find $\bmu_h^{n+1,k} \in V_h$ such that for all $v_h \in V_h$
\begin{align}
\label{eq:fixedstress2} 2G \left ( \boldsymbol \e(\bmu_h^{n+1,k}),   \boldsymbol \e(\bmv_h) - \boldsymbol \e(\bmu_h^{n+1,k} )   \right ) + \left ( \lambda \div \bmu_h^{n+1,k},   \div \bmv_h - \div \bmu_h^{n+1,k}  \right ) \\
\nonumber - \alpha \left ( {p_h^{n+1,k}},   \div \bmv_h - \div \bmu_h^{n+1,k} \right )  \geq  \left (\boldsymbol f_0,  \bmv_h - \bmu_h^{n+1,k} \right) +  \left (\boldsymbol f_2,  \bmv_h - \bmu_h^{n+1,k} \right)_{\Gamma_2}.
\end{align}
\textcolor{black}{In the flow step \eqref{eq:nf_eq4_fs_sr}, we add the regularization term $\frac{\alpha^2}{\lambda} \big(\frac{p_{h}^{n+1,k} - p_{h}^{n,k}}{\Delta t}\big)$ which is consistent. The first part $\frac{\alpha^2}{\lambda} \big(\frac{p_{h}^{n+1,k}}{\Delta t}\big)$ is added to the left hand side and a similar term on the right hand side except it is iteration lagged, that is, $\frac{\alpha^2}{\lambda} \big(\frac{p_{h}^{n,k}}{\Delta t}\big)$. In the limit (when the iteration index approaches infinity: $n \rightarrow \infty$), the terms on the left and the right hand sides will cancel, retaining the consistency of the scheme.
\noindent Our main result is the contraction of the above scheme. 
\begin{theorem} \label{thm:nf_contraction_fs_sr}
Let $\beta = \frac{1}{M \alpha^2 } + \frac{c_f}{\alpha^2} \varphi_0 + \frac{1}{\lambda}$, the iterative scheme defined by \eqref{eq:nf_eq4_fs_sr} - \eqref{eq:fixedstress2} is a contraction given by
\begin{align*}
\begin{array}{lll}
\Big\| \delta {\sigma}_{v}^{n+1,k} \Big\|_{\Omega}^{2}  + 
 \frac{2 \Delta t }{\mu_f \beta} \Big\| \boldsymbol K^{-1/2} \delta {\z}_{h}^{n+1,k}  \Big\|_{\Omega}^{2} + 4 G \lambda  \big \| \boldsymbol \e( \delta {\bmu}_{h}^{n+1,k}) \big \|_{\Omega}^{2} \\ + {\lambda}^2 \big \|\nabla \cdot \delta {\bmu}_{h}^{n+1,k} \big\|_{\Omega}^{2}
 \leq  \Big(\frac{1}{\lambda \beta}\Big)^{2} \Big\| \delta {\sigma}_{v}^{n,k} \Big\|_{\Omega}^{2}.
\end{array}
\end{align*}
\end{theorem}
}

\begin{remark}
   As stated above, the splitting scheme regularizes flow equation with the factor $\dfrac{\alpha^2}{\lambda}$ multiplied to $\frac{p_{h}^{n+1,k} - p_{h}^{n,k}}{\Delta t}$. This particular value of regularization parameter is often used in practice, however, as shown in \cite{MikelicWheeler2012} for linear Biot model, this value is not the only one. In fact, one can take any value larger than $\alpha^2/2\lambda$ to ensure convergence. A smaller value of this parameter is expected to have smaller contraction coefficient and hence faster convergence. But this is not always the case in practice. The choice of optimal value has been investigated in \cite{storvik2019optimization}.
\end{remark}
\vspace{-1.0cm}
\subsection{Proof of Contraction} 
Recall that $\beta$ is defined as: $\beta = \frac{1}{M \alpha^2 } + \frac{c_f}{\alpha^2} \varphi_0 + \frac{1}{\lambda},$ which represents the coefficient in front of the first term on the left hand side of flow equation. The volumetric mean total stress can be defined as follows \textcolor{black} {(we refer to \cite{MikelicWheeler2013} and equation (3.11) in \cite{kim2010_PhD} for the full derivation of the volumetric mean total stress)}:
\begin{align}
{\sigma}_{v}^{n,k} = {\sigma}_{v}^{k-1} + \lambda \nabla \cdot {\bmu}_{h}^{n,k} - \alpha ({p}_{h}^{n,k}-{p}_{h}^{k-1}). 
\label{eq:nf_volStressN_disc_fs_sr} 
\end{align}
where $k$ indexes time-step iteration, and $n$ indexes flow and geomechanics coupling iteration. \textcolor{black} {In terms of differences of two consecutive iterates, this can be written as:
\begin{align}
\nonumber \delta {\sigma}_{v}^{n,k} &= {\sigma}_{v}^{n,k} - {\sigma}_{v}^{n-1,k}= \lambda \nabla \cdot ({\bmu}_{h}^{n,k} - {\bmu}_{h}^{n-1,k}) - \alpha ({p}_{h}^{n,k} - {p}_{h}^{n-1,k})  \\
& =  \lambda \nabla \cdot \delta {\bmu}_{h}^{n,k} - \alpha \delta {p}_{h}^{n,k}. 
\label{eq:nf_volStressDiff1_disc_fs_sr} 
\end{align} 
}
\begin{itemize}
\vspace{-0.5cm}
\item {\bf{Step (1): Flow equations}} \\
\textcolor{black} {Consider the pressure equation \eqref{eq:nf_eq4_fs_sr}  for two successive iterates ${p}_{h}^{n+1,k}$ and ${p}_{h}^{n,k}$. Using the linearity of flow equations, we subtract the two equations to obtain equations in terms of $\delta {p}_{h}^{n,k}$. We also note that the equation in terms of $\delta {p}_{h}^{n,k}$ is linear and of the same form as the original pressure equation with no contribution of the source term $(q_h, \theta_h)$, assuming that the source term $q_h$ is fixed during iterative coupling iterations. } Now choose ${\theta}_h = \delta {p}_{h}^{n+1,k}$ to obtain, 
\begin{align*}
\Big(\frac{1}{\Delta t} {\alpha}^2 \beta \delta {p}_{h}^{n+1,k}, \delta {p}_{h}^{n+1,k} \Big) + \frac{1}{\mu_f} \Big( \nabla \cdot \delta \z_{h}^{n+1,k}, \delta p_{h}^{n+1,k}\Big) = -\frac{\alpha}{\lambda \Delta t} \Big(  \delta \sigma_{v}^{n,k}, \delta p_{h}^{n+1,k}\Big)
\end{align*}
leading to
\begin{align*}
\frac{\beta}{\Delta t} \Big\| \alpha \delta p_{h}^{n+1,k}\Big\|^2 + \frac{1}{\mu_f} \Big( \nabla \cdot \delta \z_{h}^{n+1,k}, \delta p_{h}^{n+1,k}\Big) &= \frac{1}{\Delta t} \Big(-\alpha \varepsilon \delta p_{h}^{n+1,k}, \frac{1}{\varepsilon \lambda} \delta \sigma_{v}^{n,k} \Big) \\
\nonumber \hspace{-2cm}  \leq \frac{1}{\Delta t} \Big( \frac{{\varepsilon}^2}{2} \Big\| \alpha \delta p_{h}^{n+1,k}\Big\|^2 + \frac{1}{2 {\varepsilon}^2 {\lambda}^2} \Big\|   \delta {\sigma}_{v}^{n,k}  \Big\|^2 \Big).
\end{align*}
Letting ${\varepsilon}^2 = \beta$, we obtain
\begin{align}
\label{eq:nf_eq10_fs_sr}
\beta \Big\| \alpha \delta p_{h}^{n+1,k}\Big\|^2  
+ \frac{2 \Delta t}{\mu_f}\big(\nabla \cdot \delta {\z}_{h}^{n+1,k},\delta p_{h}^{n+1,k} \big) &\leq  \frac{1}{ \beta {\lambda}^2} \Big\| \delta {\sigma}_{v}^{n,k} \Big\|^2.
\end{align}

Now, we consider the flux equation \eqref{eq:nf_eq5_fs_sr}. \textcolor{black}{As in the pressure equation case, we repeat the same process of writing the equation for two successive iterates and then subtract to obtain the equation in terms of $ \delta {\z}_{h}^{n+1,k}$. Similar to the pressure equation case, the gravity term is independent of the iterative coupling iteration $n$ and gets cancelled.} Now, we choose ${\boldsymbol q}_h = \delta {\z}_{h}^{n+1,k}$ to get
\begin{equation}
 \big(\boldsymbol K^{-1}\delta {\z}_{h}^{n+1,k},\delta {\z}_{h}^{n+1,k} \big) = ( \delta p_{h}^{n+1,k}, \nabla \cdot \delta {\z}_{h}^{n+1,k}).
\label{eq:nf_eq11_fs_sr}
\end{equation}
Substituting \eqref{eq:nf_eq11_fs_sr} into \eqref{eq:nf_eq10_fs_sr}, we obtain
\begin{align}
\label{eq:nf_eq12_fs_sr}
\beta \Big\| \alpha \delta p_{h}^{n+1,k}\Big\|^2  + \frac{2 \Delta t}{\mu_f} \Big\| \boldsymbol K^{-1/2} \delta {\z}_{h}^{n+1,k}  \Big\|^{2} &\leq  \frac{1}{ \beta {\lambda}^2} \Big\| \delta {\sigma}_{v}^{n,k} \Big\|^2.
\end{align}

\item {\bf{Step (2): Elasticity equation}}\\
Considering the variational inequality for iterative coupling iteration $n+1$ and choosing $\bmv_h = \bmu_h^{n,k}$ as the test function,  we obtain
\begin{align*}
 2G \left ( \boldsymbol \e(\bmu_h^{n+1,k}),  \boldsymbol \e(\bmu_h^{n,k}) - \boldsymbol \e(\bmu_h^{n+1,k} )   \right ) + \left ( \lambda \div \bmu_h^{n+1,k},   \div \bmu_h^{n,k} - \div \bmu_h^{n+1,k}  \right ) \\
\nonumber - \alpha \left ( p_h^{n+1,k},   \div \bmu_h^{n,k} - \div \bmu_h^{n+1,k} \right )  \geq  \left (\boldsymbol f_0,  \bmu_h^{n,k} - \bmu_h^{n+1,k} \right) +  \left (\boldsymbol f_2,  \bmu_h^{n,k} - \bmu_h^{n+1,k} \right)_{\Gamma_2}.
\end{align*}
Writing the same for $n$ and choosing the test function $\bmv_h = \bmu_h^{n+1,k}$ we get
\begin{align*}
 2G \left ( \boldsymbol \e(\bmu_h^{n,k}),   \boldsymbol \e(\bmu_h^{n+1,k}) - \boldsymbol \e(\bmu_h^{n,k} )   \right ) + \left ( \lambda \div \bmu_h^{n,k},   \div \bmu_h^{n+1,k} - \div \bmu_h^{n,k}  \right ) \\
\nonumber - \alpha \left ( p_h^{n,k},   \div \bmu_h^{n+1,k} - \div \bmu_h^{n,k} \right )  \geq  \left (\boldsymbol f_0,  \bmu_h^{n+1,k} - \bmu_h^{n,k} \right) +  \left (\boldsymbol f_2,  \bmu_h^{n+1,k} - \bmu_h^{n,k} \right)_{\Gamma_2}.
\end{align*}
Adding the two equations together we write using our notation for the difference of two consecutive iterations,

\begin{align}
\label{eq:nf_eq13_fs_sr}
2 G   \big \| \boldsymbol \e( \delta {\bmu}_{h}^{n+1,k})\big \|^2 + \lambda  \big \|\nabla \cdot \delta {\bmu}_{h}^{n+1,k} \big\|^2 
- \alpha \big( \delta p_{h}^{n+1,k} ,\nabla \cdot \delta {\bmu}_{h}^{n+1,k} \big) \leq 0. 
\end{align}
Multiplying both sides of \eqref{eq:nf_eq13_fs_sr} by \textcolor{black}{$(2 \lambda \beta)$}, we get
\textcolor{black} {
\begin{align}
\label{eq:nf_eq14_fs_sr} &4 G \lambda \beta  \big \| \boldsymbol \e( \delta {\bmu}_{h}^{n+1,k}) \big \|^2 + 2\lambda^2 \beta \big \|\nabla \cdot \delta {\bmu}_{h}^{n+1,k} \big\|^2 \\
\nonumber &\hspace{2.0cm} \leq 2\alpha \lambda \beta \big( \delta p_{h}^{n+1,k} ,\nabla \cdot \delta {\bmu}_{h}^{n+1,k} \big).
\end{align}
Multiplying \eqref{eq:nf_eq12_fs_sr} and \eqref{eq:nf_eq14_fs_sr} by $\frac{1}{\beta}$  and adding them up, we get
\begin{align*}
&\Big\| \alpha \delta p_{h}^{n+1,k}\Big\|^2  + \frac{2 \Delta t }{\mu_f \beta} \Big\| \boldsymbol K^{-1/2} \delta {\z}_{h}^{n+1,k}  \Big\|^{2} + 4 G \lambda  \big \| \boldsymbol \e( \delta {\bmu}_{h}^{n+1,k}) \big \|^2 \\  
&+ 2 {\lambda}^2 \big \|\nabla \cdot \delta {\bmu}_{h}^{n+1,k} \big\|^2 - 2 \big(\alpha \delta p_{h}^{n+1,k} ,\lambda \nabla \cdot \delta {\bmu}_{h}^{n+1,k} \big) \leq  \Big(\frac{1}{\lambda \beta}\Big)^{2} \Big\| \delta {\sigma}_{v}^{n,k} \Big\|^2
\end{align*}
which can be written as
\begin{align*}
&\Big\{ \Big\| \alpha \delta p_{h}^{n+1,k}\Big\|^2 - 2 \big(\alpha \delta p_{h}^{n+1,k} ,\lambda \nabla \cdot \delta {\bmu}_{h}^{n+1,k} \big)  +  \big \| \lambda \nabla \cdot \delta {\bmu}_{h}^{n+1,k} \big\|^2\Big\} \\  
& + \frac{2 \Delta t }{\mu_f \beta} \Big\| \boldsymbol K^{-1/2} \delta {\z}_{h}^{n+1,k}  \Big\|^{2} + 4 G \lambda  \big \| \boldsymbol \e( \delta {\bmu}_{h}^{n+1,k}) \big \|^2 \\
& \hspace{0.3cm}+ {\lambda}^2 \big \|\nabla \cdot \delta {\bmu}_{h}^{n+1,k} \big\|^2
\leq  \Big(\frac{1}{\lambda \beta}\Big)^{2}\Big\| \delta {\sigma}_{v}^{n,k} \Big\|^2.
\end{align*}
Thus we have
\begin{align}
\nonumber &  \Big\| \delta {\sigma}_{v}^{n+1,k} \Big\|^2 + \frac{2 \Delta t }{\mu_f \beta} \Big\|\boldsymbol K^{-1/2} \delta {\z}_{h}^{n+1,k}  \Big\|^{2} + 4 G \lambda  \big \| \boldsymbol \e( \delta {\bmu}_{h}^{n+1,k}) \big \|^2 \\
\label{eq:nf_eq16_fs_sr} & \hspace{0.3cm} + {\lambda}^2 \big \|\nabla \cdot \delta {\bmu}_{h}^{n+1,k} \big\|^2
\leq  \Big(\frac{1}{\lambda \beta}\Big)^{2} \Big\| \delta {\sigma}_{v}^{n,k} \Big\|^2.
\end{align}
}
\end{itemize}
 It is clear that the contraction constant: $\Big(\frac{1}{\lambda \beta}\Big)^{2} = \Big(\frac{M {\alpha}^2}{\lambda + M \lambda c_f \varphi_0 + M {\alpha}^2}\Big)^{2} < 1$.
This completes the proof. 

\begin{remark}
Note that our contraction result is in terms of the composite quantities that are a combination of $\bmu_h^{n,k}$ and $p_h^{n,k}$. We need to ensure that this leads to the convergence of  the unknown quantities in our discrete form. This is indeed ensured by the existence of  limit functions ${p}_{h}^{k},  {\bmu}_{h}^{k}, {\z}_{h}^{k}$ such that
\begin{align*}
\begin{array}{llllllllllllllll}
 {p}_{h}^{n,k} \rightarrow {p}_{h}^{k}  & \text{ in } L^2({\Omega}),  &   {\bmu}_{h}^{n,k}   \rightarrow {\bmu}_{h}^{k} & \text{ in } {H^1({\Omega})}^d , & {\z}_{h}^{n,k} \rightarrow   {\z}_{h}^{k}  & \text{ in } {H(\rm{div},\Omega})^d 
\end{array}
\end{align*}
converge strongly in the norms of the above spaces.
\end{remark}

We outline the argument for the above assertion. 
The contraction result in \eqref{eq:nf_eq16_fs_sr} implies that $\big \|  \delta \sigma_{v}^{n+1,k}  \big \|_{\Omega}$, $\big \|\nabla \cdot  \delta {\bmu}_{h}^{n+1,k} \big\|_{\Omega}$, and $\big \| \boldsymbol K^{-1/2}\delta  {\z}_{h}^{n+1,k} \big\|_{\Omega}$ converge geometrically to zero. This implies that ${\sigma}_{v}^{n+1,k}$, $\nabla \cdot {\bmu}_{h}^{n+1,k}$, and ${\z}_{h}^{n+1,k}$ are Cauchy sequences in $L^2(\Omega)$. By the definition of $\sigma_v^{n,k}$, i.e. \eqref{eq:nf_volStressN_disc_fs_sr} and \eqref{eq:nf_volStressDiff1_disc_fs_sr}, we obtain that $p_h^{n,k}$ is a Cauchy sequence in Hilbert space and has a unique limit in $L^2({\Omega})$. Similarly, for the displacement, \eqref{eq:nf_eq16_fs_sr} implies that $\e( \delta {\bmu}_{h}^{n+1,k})$ converges geometrically to 0 in $L^2({\Omega})$, which implies that ${\bmu}_{h}^{n+1,k}$ converge geometrically in ${H^1({\Omega})}^d$. Thus, ${\bmu}_{h}^{n+1,k}$ is a Cauchy sequence in a complete Hilbert space, and hence has a unique limit in the corresponding space. For the divergence of the flux, we note 
 both $\nabla \cdot {\z}_{h}^{n+1,k}$ and  ${\z}_{h}^{n+1,k}$ converge geometrically to zero in $L^2({\Omega})$, hence ${\z}_{h}^{n+1,k}$ converges in ${H(\rm{div},\Omega})^d$. The existence of the limit follows from the completeness of the space. 
\begin{remark}
    Even though we carried the convergence analysis for a particular spatial discretization scheme, the above convergence proof can be easily adapted to other schemes. This is evident from the fact that the proof does not use any particular properties of discretization. For example, a similar analysis of conformal Galerkin discretization for both flow and mechanical deformation immediately follows.  
\end{remark}
\section{Conclusion and future work}
 \textcolor{black} {In this work, we developed and analysed the convergence of a fixed stress split scheme for a linear poroelastic problem with contact conditions. The scheme was shown to be contractive. For future work, we plan to include the friction behaviour into the model and run numerical simulations to study the performance of the fixed stress scheme. }


\begin{acknowledgement}
KK acknowledges the support of the VISTA program, The Norwegian Academy of Science and Letters, and Equinor.    
\end{acknowledgement}


\bibliographystyle{unsrt}
\bibliography{Fixedstress_contact_Biot}

\end{document}